\documentclass[12pt]{amsart}
\newfont{\sheaf}{eusm10 scaled\magstep1}
\usepackage{amssymb} 
\usepackage[mathscr]{eucal} 
\usepackage[all]{xy}
\usepackage{amsmath} 

\newcommand{\C}{\ensuremath{\mathbb{C}}}
\newcommand{\R}{\ensuremath{\mathbb{R}}}
\newcommand{\Z}{\ensuremath{\mathbb{Z}}}

\newcommand{\hol}{\ensuremath{\mathcal{O}}}

\newcommand{\PP}{\ensuremath{\mathbb{P}}}
\newcommand{\RR}{\ensuremath{\mathcal{R}}}

\newcommand{\FF}{\ensuremath{\mathcal{F}}}
\newcommand{\HHH}{\ensuremath{\mathcal{H}}}

\newcommand{\Proof}{{\it Proof. }}
\newcommand{\ra}{\ensuremath{\rightarrow}}

\newcommand{\sC}{{\mathcal C}}
\newcommand{\sT}{{\mathcal T}}

\newcommand{\sM}{{\mathcal M}}

\newcommand{\la}{\lambda}
\newcommand{\ze}{\zeta}
\newcommand{\Ga}{\Gamma}

\newcommand{\La}{\Lambda}
\newcommand{\ga}{\gamma}
\newcommand{\e}{\epsilon}
\newcommand{\into}{\hookrightarrow}

\newcommand{\QED}{\hspace*{\fill}$Q.E.D.$}

\def\eea{\end{eqnarray*}}
\def\bea{\begin{eqnarray*}}

\def\hol{{\mathcal{O}}}

\newtheorem{teo}{Theorem}[section]
\newtheorem{df}[teo]{Definition}
\newtheorem{lem}[teo]{Lemma}
\newtheorem{cor}[teo]{Corollary}
\newtheorem{oss}[teo]{Remark}
\newtheorem{prop}[teo]{Proposition}

\newtheorem{REM}[teo]{Remark}

\def\eea{\end{eqnarray*}}
\def\bea{\begin{eqnarray*}}

\begin{document}

\title{Deformation in the large of some complex manifolds, II}

\author{Fabrizio Catanese\\
    Paola Frediani\\
    }

\footnote{
The present research took place in the framework of the Schwerpunkt
"Globale Methode in der komplexen Geometrie", and of the PRIN 2003:
"Spazi dei moduli e
    teoria di Lie" (MURST)

 }
\date{May 18, 2005}
\maketitle

\begin{abstract}
The compact complex manifolds considered in this article are
principal torus bundles over a torus. We consider the
Kodaira Spencer map of the complete Appell Humbert family
(introduced by the first author in Part I) and are able to show
that we obtain in this way a connected component of the space
of complex structures each time that the base dimension is two, the fibre
dimension is one, and a suitable topological condition is verified.
\end{abstract}

\section{Introduction}

    One of the most interesting problems in the theory of compact
complex manifolds
is " Moduli theory", i.e., the study of the  space of complex
structures $\sC (M)$
on a given oriented differentiable manifold $M$. Viewing these as integrable
almost complex structures yields a description as an infinite dimensional
space on which  acts the infinite dimensional group $ Diff^+(M)$ of
orientation preserving diffeomorphisms of $M$.

Moduli theory in a proper sense means studying in detail the
quotient $\frak M (M) : =  \sC (M) / Diff^+(M) $, and often this
study is split into
two parts.

First one studies the Teichm\"uller space of $M$, i.e., $ \frak T
(M): =  \sC (M) /
Diff^0(M)
$, where $ Diff^0(M) \subset Diff^+(M) $ is the connected component of the
identity (its elements are called the diffeomorphisms {\bf isotopic}
to the identity). Then one may view the Moduli space  $\frak M (M)$ as
the quotient of the Teichm\"uller space $ \frak T (M)$ by the Mapping class
group $ \sM ap (M) : = Diff^+(M) / Diff^0(M)$, hoping that this action
turn out to be properly discontinuous.

Of course a preliminary question is   the determination of the connected
components of $\sC (M)$, which are called the
deformation classes
in the large of  the complex structures on $M$, and which we believe may be
countably many (their number may  be arbitrarily high, cf.
\cite{cat02}, \cite{cat04}).
Observe that
the usual deformation theory addresses only the study of the {\bf small
deformations}, i.e., it describes the so called Kuranishi space $\frak B (M)$
which is the germ of
the  Teichm\"uller space at the point corresponding to a given
complex structure.
Kuranishi's theorem (\cite{ku1},\cite{ku2}) shows that the
Teichm\"uller space is a
complex space  of locally finite dimension
(its dimension may however be unbounded, cf.
cor. 7.7 of
\cite{cat04}).

A prototype test for all these desiderata  is furnished by the example of
the complex tori.
These are parametrized  by an open set $\mathcal T_n$ of the complex
Grassmann Manifold $Gr(n,2n)$, image of the open set of matrices

$\{ \Omega \in Mat(2n,n; \C) \ | \ (i)^n det  (\Omega \overline 
{\Omega}) > 0 \}.$

This parametrization is very explicit: if we consider
a fixed lattice $ \Ga \cong \Z^{2n}$,
to each matrix $ \Omega $ as above we associate the subspace $ V = ( \Omega )
(\C^{n})$, so that
$ V \in Gr(n,2n)$ and $\Ga \otimes \C \cong V \oplus \bar{V}.$

Finally, to $ \Omega $ we associate the torus $Y_V : = V / p_V (\Ga)$,
$p_V : V \oplus
\bar{V} \ra V$ being the projection onto the first addendum.

Not only we obtain in this way a connected open set
inducing  all the small deformations (cf.
\cite{k-m71}), but indeed, as  it was shown in \cite{cat02} (cf. also
\cite{cat04})
a connected component of the Teichm\"uller space.

It was observed however by Kodaira and Spencer already in their first article
(\cite {k-s58}, and volume II of Kodaira's
collected works)
that the mapping class group $ SL ( 2n, \Z)  $ does not act properly
discontinuously
on $\mathcal T_n$. This shows that for compact complex manifolds
it is better to consider, rather than the Moduli space, the
Teichm\"uller space,
which can be obtained by glueing together several  Kuranishi spaces.

Moreover, after some initial constructions
by Blanchard and Calabi (cf. \cite{bla53}, \cite{bla56},
\cite{cal58}) of non K\"ahler
complex structures
$X$ on manifolds diffeomorphic to a product $C \times T$, where $C$ is a
compact complex curve and $T$ is a 2-dimensional complex torus,
Sommese observed (\cite{somm75}) that  the space of complex
structures on a  six dimensional real torus is not connected.

These examples were then generalized in
\cite{cat02} \cite{cat04} under the name of {\bf Blanchard-Calabi manifolds}
showing (corollary 7.8 of \cite{cat04}) that also
the space of complex structures on the product of a curve $C$
of genus $g \geq 2$ with a four
dimensional real torus is not connected.

On the positive side, however, our goal is to determine completely
the connected components of spaces $\sC (M)$ (the so called equivalence
classes for "deformation in the large") for
suitable differentiable manifolds.

    A very interesting class of examples to consider is given by the quotients
of an affine
space  $\C^n$ by a nilpotent
or solvable group $\Pi$.
One of the main tools in order to analyse
this  type of examples should be the cited result that every
deformation in the large of a complex torus is a complex torus.

Since this ambitious program carries several big difficulties with it,
we begin to address it in the simplest possible case, i.e., when
we have a holomorphic torus bundle
$f :  X \ra Y$ over a torus $Y$ and the group
is 2-step nilpotent (the central group extension is given by
$  \Lambda : = \pi_1(F) \into \pi_1(X) = \Pi \ra \pi_1(Y) =
\Gamma$, classified by an alternating bilinear form
$ A \in \Lambda^2(\Gamma  )^{\vee} \otimes (\Lambda)$).

We shall use the classification theory for principal
holomorphic torus bundles  over tori  given in Theorem 6.8 of \cite{cat04},
where  an explicit family is produced, described purely in terms
of bilinear algebra, and called the {\bf complete Appell-Humbert family
of Torus Bundles}.

Roughly speaking,
one looks for the subspaces $U \subset \La \otimes \C,
V \subset \Ga \otimes \C$ such that,
viewing $A$ as a real element of
    $$\Lambda^2(\Gamma
\otimes \R )^{\vee} \otimes (\Lambda \otimes \R )=  \Lambda^2(V \oplus
\bar{V})^{\vee} \otimes (U \oplus \bar{U}),$$
    its component in $  \Lambda^2(
\bar{V})^{\vee} \otimes (U )$ is zero (one says then that the
{\bf Riemann bilinear relations} hold).

This theory bears close similarities with the theory of line bundles
over tori, and allows for explicit descriptions of several holomorphic
invariants (see theorems 6.11-6.12 of \cite{cat04}).

The simplest way to explain the meaning of the Appell- Humbert
family is as follows: one observes that the differentiable manifold
$M$ underlying $X$ is
simply the quotient of a contractible real Lie group $\Pi \otimes \R$ by a
discrete subgroup
$\Pi $.  We consider now all the right invariant almost complex structures
on $\Pi \otimes \R$, and we show that the Riemann bilinear relations are
equivalent to the integrability of these almost complex structures.

    In particular,  explicit formulae can be given for this family of
complex structures, and from these one derives for instance

{\bf  Theorem 6.10 of \cite{cat04} :}
{ \it
The cokernel of $ 0 \ra H^0 (\Omega^1_Y )\ra H^0 (\Omega^1_X )$ is the
subspace of ${U}^{\vee}$ which annihilates the image of the Hermitian
part $B''$  of $A$, $B'' \in [( V \otimes \bar{V})^{\vee}  \otimes (U
)]$, i.e.,
$\{ {u}^{\vee} \in {U}^{\vee} | {u}^{\vee} \dashv B'' = 0 \in ( V \otimes
\bar{V})^{\vee} \}$.

It follows in particular that $X$ is parallelizable if and only if the
Hermitian part of $A$ is zero.}

We already observed that $X$, as a differentiable manifold, is the
quotient of the real Lie group $\Pi \otimes \R$  by the discrete
subgroup $\Pi $: the point is that  there are some complex structures
where we can realize $X$ as the quotient of a complex Lie group by
$\Pi $.

The most famous example is the Iwasawa 3-fold (cf. \cite{k-m71},
whose small deformations were analysed
     by Nakamura in
\cite{nak75}, who thus showed that small deformations of a complex
parallelizable
manifold need not be parallelizable.

Our first result in this paper will be the following

{\bf Theorem 4.11}
{\em Let $ f:  X \ra Y$ be a holomorphic principal bundle with base a
complex torus
$ Y: = Y _V$ of dimension $m = 2$, and fibre an elliptic curve $T : = T_U$.

Then the complete Appell Humbert family is versal at $X$ if
either $B' \neq 0$, or $B'=0$ and  $B (v, x)$ non degenerate in the
first factor,
or $X$ is parallelizable and $A \neq 0$. }

The strong restriction on the dimensions $m=2, d=1$ is forced
by the requirement that the complete Appell Humbert family have
a smooth basis. This condition,  combined with the
surjectivity of the  Kodaira Spencer map, ensures that all the small
deformations
are again principal holomorphic torus bundles.

The above result combines with the following

{\bf Theorem 6.17 of \cite{cat04}}
{\em Let $ f:  X \ra Y$ be a holomorphic principal bundle with base a
complex torus
$ Y _V$ of dimension $m$, and fibre an elliptic curve $T_U$.

Assume moreover $ \pi_1 (X) := \Pi $ to be a nontrivial  central extension
$$ 1 \ra \Lambda \ra \Pi \ra \Gamma \ra 1$$
classified
by  a cohomology class
$\epsilon  \neq 0 \ \in H^2(Y, \Lambda )$
whose associated bilinear form $A$ has an image of dimension $=2$.

Then every limit of  manifolds in the complete Appell-Humbert family is again a
holomorphic principal bundle $ f' :  X' \ra Y'$ with fibre an
elliptic curve $T'$,
and thus occurs in the complete Appell-Humbert family .}

    and one shows that the class
of holomorphic torus bundles with
    fibre dimension
$d=1$,  and base dimension $m=2$ form irreducible components
in the Teichm\"uller space.

With a further assumption, however, we are able to produce many
connected components
of the Teichm\"uller space

{\bf Main Theorem 4.13}
{\em Consider the family of holomorphic principal torus bundles $ f:  X \ra Y$
with  base a complex torus
$ Y _V$ of dimension $2$, and with  fibre an elliptic curve $T_U$,
corresponding to a bilinear form $A$ such that

1) $A$ is non degenerate and ${\rm Im}  A$ has dimension $2$

2) the associated pencil

$ \PP^1 \cong \PP ( ( \La \otimes \R)^{\vee} ) \into \PP
   ( \La^2(\Ga \otimes
\R)^{\vee}) \cong
\PP ( \La^2( \R^{4}))$
intersects the Klein Pfaffian quadric in at least one real point.

Then this family forms a connected component of the  Teichm\"uller space.}

Indeed condition 2) is necessary: in fact there are deformations of
the Iwasawa manifold for which the complete Appell Humbert family
has not surjective Kodaira Spencer map.

Time reasons do not allow us to investigate fully the difficult
question whether
in this and similar cases the complete Appell Humbert family fails to
yield an open set
in the  Teichm\"uller space.

\section{Principal holomorphic torus bundles over tori: generalities}

Throughout  the paper, our set up will be the
following: we have a holomorphic submersion between compact
complex manifolds

$$ f : X \ra Y  ,$$
such that the base $Y$ is a complex torus, and one fibre $F$  (whence
all the fibres, by
theorem 2.1 of
\cite{cat04}) is also a  complex torus.

We shall denote this situation by saying that $f$ is
differentiably a torus bundle.

We let $n= dim X$, $m = dim Y$, $d=dim F= n -m$.

In general  ( cf.\cite {fg65})) $f$ is a
holomorphic bundle if and only if all the  fibres are
biholomorphic.

This holds necessarily in the special case $d=1$, because the moduli space
for $1$-dimensional complex tori exists and is isomorphic to $\C$.

Assume now that we have a holomorphic torus fibre bundle, thus we have
    (cf. \cite{cat04}, pages 271-273) the exact sequence
$$ \  0 \ra  \Omega^1_Y \ra f_* \Omega^1_X \ra f_* \Omega^1_{X|Y} \ra 0 .$$
and a principal
holomorphic bundle if moreover $f_* \Omega^1_{X|Y}$ is holomorphically trivial.

\begin{oss}
In general (cf. e.g. \cite{bpv84}) if $T$ is a complex torus, we have an
exact sequence of complex Lie groups
$$ 0 \ra T \ra Aut(T) \ra M \ra 1 $$
where $M$ is discrete.
Taking sheaves of germs of holomorphic maps with source $Y$ we get
$$ 0 \ra \HHH(T)_Y \ra \HHH (Aut(T))_Y \ra M \ra 1 $$
and the exact sequence
$$ 0 \ra H^1(Y, \HHH (T)_Y) \ra H^1(Y, \HHH (Aut(T)))_Y \ra H^1(Y,
M) $$
(since holomorphic bundles with base $Y$ and fibre $T$ are
classified by the cohomology group $H^1(Y, \HHH (Aut(T))_Y)$)
determines the discrete obstruction for a holomorphic bundle
    to be a principal holomorphic bundle.

\end{oss}

In the case of a principal holomorphic bundle we
write $\pi_1(T ) \cong \Lambda$,  $\pi_1(Y ) \cong \Ga$
and the exact sequence
$$\   \ra H^0(\HHH(T)_Y) \ra H^1(Y, \Lambda ) \ra  H^1 (Y,
\hol_Y^d) \ra H^1(\HHH(T)_Y) \ra ^c \ra H^2(Y, \Lambda ) $$
determines a  cohomology class
$\epsilon \in H^2(Y, \Lambda )$
which classifies the  central extension
$$  1 \ra \pi_1(T) = \La \ra \Pi := \pi_1 (X) \ra \pi_1(Y) = \Ga \ra 1$$
(it is central by the triviality of the monodromy automorphism).

\begin{prop}

\label{unicover}

Let $f: X \rightarrow Y$ be a principal holomorphic torus bundle over
a torus as
above.

Then the universal covering of $X$ is isomorphic to ${\C}^{m+d}$ and $X$ is
biholomorphic to a quotient $X \cong {\C}^{m+d}/\Pi$.

\end{prop}

{\bf Proof.}

Let $\pi: {\C}^m \rightarrow Y$ be the universal covering map and let
us consider
the pull back bundle ${\pi}^*X \rightarrow {\C}^m$ as in the
following Cartesian
diagram:

$$\xymatrix{
{\pi}^* X \ar[r]^{{\pi'}} \ar[d]&X  \ar[d]\\
{\C}^m  \ar[r]^{{\pi}}&  Y }$$

$\pi^*X \rightarrow {\C}^m$ is a principal holomorphic torus bundle on
${\C}^m$, therefore it is trivial.  In fact principal holomorphic
torus bundles on
${\C}^m$ are classified by $H^1({\C}^m , {\mathcal
H}(T)_{{\C}^m})$,and we claim
that
$H^1({\C}^m , {\mathcal H}(T)_{{\C}^m}) =0 $. Let us in fact consider the exact
sequence defining $T$,

$$1 \rightarrow \Lambda \rightarrow U \rightarrow T \rightarrow 1$$

Taking germs of holomorphic functions with source ${\C}^m$, we get

$$1  \rightarrow \Lambda \rightarrow {\mathcal H}(U)_{{\C}^m} \rightarrow
{\mathcal H}(T)_{{\C}^m} \rightarrow 1$$

and the long exact cohomology sequence  yields

$$0=H^1({\C}^m, \Lambda) \rightarrow H^1({\C}^m, {\mathcal
O}^d_{{\C}^m})=0 \rightarrow H^1({\C}^m,  {\mathcal H}(T)_{{\C}^m})
\rightarrow H^2({\C}^m, \Lambda)=0,$$

therefore $H^1({\C}^m, {\mathcal H}(T)_{{\C}^m}) =0$.

So we have $\pi^*X \cong {\C}^m \times T$ and by taking the universal covering
${\C}^m \times {\C}^d \rightarrow \pi^*X \cong {\C}^m \times T$  and the
composition with the covering $\pi' : \pi^*X \rightarrow X$, we
obtain the desired
assertion.
\hfill{Q.E.D.}

Let us briefly recall again the standard  way to look at the family
$\sT_m$ of complex
tori of complex dimension $=m$. We fix a lattice $\Ga$ of rank $ 2m$, and
we look at the complex ($m$-dimensional) subspaces $ V \subset \Ga \otimes \C$
such that  $ V \oplus \bar{V} = \Ga \otimes \C$: to $V$ corresponds the
complex torus $ Y_V : = \Ga \otimes \C / (\Ga \oplus \bar{V} )$.
We finallyselect one of the two resulting connected components by requiring
that the complex orientation of $V$ induces on $\Ga \cong p_V (\Ga)$
a standard orientation.

We define similarly $T_U  : = \La \otimes \C / (\La \oplus \bar{U} )$.

Consider now  our principal holomorphic torus bundle $f : X \rightarrow Y$
over a
complex torus $Y_V$ of dimension $m$, and with fibre a complex torus
    $T_U$ of dimension $d$  and let  $\epsilon \in
H^2(Y, \Lambda) = H^2(\Ga, \Lambda)$ be the cohomology class classifying
the central extension

\begin{equation}
\label{central}
1 \rightarrow \Lambda \rightarrow \Pi \rightarrow \Gamma \rightarrow 1.
\end{equation}
\begin{lem}

It is possible to "tensor" the above exact sequence with $\R$,
obtaining an exact
sequence of Lie Groups
$$  1 \ra \Lambda \otimes \R \ra \Pi \otimes \R \ra \Gamma
\otimes \R \ra 1$$ such that
$\Pi$ is a discrete subgroup of $\Pi \otimes \R$ and such that
$X$ is diffeomorphic to the quotient $$  M:=  \Pi \otimes \R /  \Pi .$$
\end{lem}

\Proof
Consider, as usual, the map
$$A: \Ga \times \Ga \ra \La,$$
$$A(\ga,\ga') =  [\hat{\ga}, \hat{\ga'}] =  \hat{\ga}
\hat{\ga'}(\hat{\ga})^{-1}(
\hat{\ga'})^{-1},$$
where $\hat{\gamma}$ and $\hat{\gamma'}$  are  respective liftings to $\Pi$ of
elements $\ga, \ga' \in \Ga$.
We observe that since the extension (\ref{central}) is central,
    the definition of $A$ does
not depend on the choice of the liftings of $\ga$, resp. $\ga'$ to $\Pi$.

As well known, $A$ is bilinear and alternating, so $A$ yields a cocycle in
$H^2(\Ga, \La)$ which "classifies" the central extension (\ref{central}).
Let us review how does this more precisely hold.

Assume  that $\{\ga_1, ..., \ga_{2m}\}$ is a basis of $\Ga$ and choose
    fixed  liftings $\hat{\ga_i}$ of $\gamma_i$ in $\Pi$, for each $i =
1, \dots , 2m$.
Then automatically we have determined a canonical way to lift
elements $\ga \in \Ga$ through:
$$ \ga = {\ga}_1^{n_1} \dots  {\ga}_{2m}^{n_{2m}} \mapsto \hat{\ga} : =
{\hat{\ga}_1}^{n_1}... {\hat{\ga}_{2m}}^{n_{2m}}.$$
Hence a canonical way to write the elements of $\Pi$ as products
    $\lambda \hat{\ga}$,
where  $\la \in \La$ and $\hat{\ga}$ is as above.

Since $\forall i,j$, one has

$$\hat{\ga_i} \hat{\ga_j} = A(\ga_i, \ga_j) \hat{\ga_j} \hat{\ga_i},$$

     we have a standard way of computing the products
$(\lambda \hat{\ga}) (\lambda' \hat{\ga}')$ as
$ \la'' (\widehat{(\ga \ga'}))$ in  $\Pi$,  where $\la ''$ will be
computed using $A$.

We can also view $\Pi$ as a group of affine transformations of $(\La
\otimes {\R})
\oplus (\Ga \otimes \R)$. In fact,   $(\La \otimes {\R})
\oplus (\Ga \otimes \R)$ is a real vector space with basis
$\{\lambda_1,..., \la_{2d},
\ga_1, ..., \ga_{2m}\}$ where $\{\lambda_1,..., \la_{2d}\}$ is a basis of
$\La$ and the action of $\Pi$ on $(\La \otimes {\R})
\oplus (\Ga \otimes \R)$
is given as follows:

$\lambda_i$ acts on $(\La \otimes {\R}) \oplus (\Ga \otimes \R)$
sending $(y, x)$ to
$(y + \lambda_i, x)$, while the action of $\hat{\gamma_j}$ is defined using the
multiplication $(\lambda \hat{\gamma})
\mapsto (\lambda \hat{\gamma}) \hat{\gamma_j}$.

More precisely if $y  \in \La \otimes \R$, $x = \sum x_j \gamma_j \in
\Ga \otimes
\R$, $\ga'  = \sum \nu_h \ga_h \in \Ga $, we have
$$(y,x)\hat{\ga'} : =  (y + \phi_{ \ga'}(x), x + \ga'),$$
where
$$\phi_{ \ga'}(x) = \sum_{j \geq h} x_j \nu_h A( \gamma_j,  \gamma_h) =
\sum_{j \geq h} x_j A_{jh} \nu_h = ^tx T^-  \ga',$$
where $T^-$ is the lower triangular part of the matrix $A$, so that
we can write
$A = T^- - ^tT^- $.

Therefore  we can endow
$(\La \otimes {\R}) \oplus (\Ga \otimes \R) =: \Pi \otimes \R$ with a Lie group
structure defined by
$$(y,x) (y', x') = (y + y' + T^-(x,x'), x + x'),$$
and the quotient $(\Pi \otimes \R)/\Pi$  of this Lie group by the
discrete subgroup $\Pi$ is easily seen to be
diffeomorphic to
$X$.

    \QED

\begin{REM}
We can change coordinates in $(\La \otimes {\R}) \oplus (\Ga \otimes
\R)$ in such a
way that the action of  the set $\hat{\Ga} \cong \Ga \subset \Ga
\otimes \R$  on $\Pi
\otimes
\R$ is given by
$$(y,x) \hat {\ga} = (y + A(x, \ga) + 2 S ({\ga}, \ga) , x + \ga),$$
where $S ({\ga}, \ga') $ is a symmetric bilinear  $(\frac{1}{4} \La )$-
valued form, and $ 2 S ({\ga}, \ga) \in \La$.
\end{REM}

\Proof

Let us define the symmetric form $S :  = - \frac{T ^-+^tT^-}{4}$, so
that $T^-  + 2S=
\frac{T^- - ^tT^-}{2} = \frac{A}{2}$.

Consider the map $ F: (\La \otimes {\R}) \oplus (\Ga \otimes \R) \ra
(\La \otimes {\R})
\oplus (\Ga \otimes \R)$ defined by $F(y,x) = (2(y + S(x,x)), x) =:( \eta, x)$.

Then $\forall \hat {\ga}  \in \hat {\Ga} $ we have an induced action
$$(\eta, x) \hat {\ga}  = F ((y,x) \hat {\ga} ) = F(y + T^-(x, \ga),
x + \ga) =$$
$$ (2y + 2 T^-(x, \ga) + 2S(x+\ga,x + \ga), x + \ga) =
(\eta + A(x, \ga) + 2S(\ga,\ga), x + \ga),$$
and we conclude observing that $2S(\ga,\ga)   = T^-(\ga,\ga) \in \La$.

\QED

We recall from \cite{cat04}  the First Riemann
bilinear Relation: it is derived from the exact cohomology
sequence
$$ H^1 (Y, \hol_Y\otimes U) \cong H^1 (Y, \HHH(U)_Y) \ra H^1(\HHH(T)_Y)
\ra ^c \ra H^2(Y, \Lambda ) \ra H^2 (Y, \HHH(U)_Y) $$
    and says that the class $\epsilon$ maps to zero in $H^2 (Y, \HHH(U)_Y) $.
More concretely we have the

{\underline{\bf First Riemann Relation for  principal holomorphic
Torus Bundles} }

{\bf Let $A : \Gamma \times \Gamma \ra \Lambda$ be the alternating bilinear map
representing the cohomology class $\epsilon$ :
    then
    $$ A \in \Lambda^2(\Gamma
\otimes \R )^{\vee} \otimes (\Lambda \otimes \R ) \subset  \Lambda^2(\Gamma
\otimes \C )^{\vee} \otimes (\Lambda \otimes \C) \subset \Lambda^2(V \oplus
\bar{V})^{\vee} \otimes (U \oplus \bar{U}),$$
satisfies the property that its  component  in $
\Lambda^2(
\bar{V})^{\vee} \otimes (U )$ is zero.}

We want now to explain in detail the bilinear algebra underlying the
Riemann bilinear
relation.

We observe preliminarly that  one has a natural isomorphism
$ \Lambda^2(V \oplus
\overline{V})^{\vee}  \cong  \Lambda^2(V)^{\vee}  \oplus
(V^{\vee} \otimes \overline{V})^{\vee}   \oplus
\Lambda^2(\overline{V})^{\vee}$,
where the middle summand embeds by the wedge product :
$ w' \otimes \bar{w} \mapsto   w' \wedge \bar{w} = w' \otimes \bar{w} - \bar{w}
\otimes w'$.

    Consider the bilinear form
$$A\in \Lambda^2(\Gamma
\otimes \C)^{\vee} \otimes (\Lambda \otimes \C ) = \Lambda^2(V \oplus
\overline{V})^{\vee} \otimes (U \oplus
\overline{U}) ,$$  satisfying the first bilinear relation and let us write
$$A = B + \overline{B},$$ where
$B \in \Lambda^2(\Gamma
\otimes \C )^{\vee}  \otimes U,$
and
$\overline{B} \in \Lambda^2(\Gamma
\otimes \C )^{\vee}  \otimes \overline{U}.$

By the first bilinear relation $B = B' + B''$, with $B' \in \Lambda^2(V)^{\vee}
\otimes U$, $B'' \in (V^{\vee} \otimes \overline{V})^{\vee} \otimes U.$

Concretely, $ A = B' + B'' + \overline { B'} +\overline { B'' }$, where
$B'$ is an alternating complex bilinear form.
The fact that $B''$ is
alternating reads out as:
$$B'' (v', \bar{v}) = - B'' (\bar{v}, v')  \   \forall v, v' \in V$$
whereas conjugation of tensors reads out as:
$$\bar{B} (\bar{x}, \bar{y}) = \overline {B (x, y)}\  \forall x,y \ \Rightarrow
\bar{B''} (\bar{v}, v') = \overline {B'' (v, \bar{v'})}.$$

\begin{df}
We define the associated  (vector-valued) Hermitian bilinear form
$$ D : V \times V \ra \La \otimes  \C$$ through
$$ D (v_1, v_2) : = i A (v_1, \bar{v_2}).$$

\end{df}

In fact $D$ is clearly complex linear in $v_1$ and complex antilinear
in $  v_2$, and
$ D (v_2, v_1) : = i A (v_2, \bar{v_1})= - i A (\bar{v_1}, v_2) = -i
\overline{A (v_1, \bar{v_2}) } = \overline{D (v_1, v_2) }.$

We have
$$ D (v_1, v_2) : = i A (v_1, \bar{v_2})= i(B''(v_1, \overline{v_2}) +
\overline{B''}(v_1, \overline{v_2})) = i B''(v_1, \overline{v_2}) +
\overline{i B''(v_2, \overline{v_1})}.$$

\begin{oss}
One can view $ D(v,v) : V \ra \La \otimes \R$ as a real  linear system
   (of dimension $ \leq 2d$) of  Hermitian forms on $V$, and its discriminant

$ \Delta_D = \{ (\la_i ) \ |  det (\sum_1^{2d}\lambda_i D_i) = 0 \} $ will be a
real hypersurface of degree $m$ in the projective space
$\PP (( \La \otimes \C)^{\vee} )$.

Observe that this linear system is independent of the
choice of $U$, but it depends upon the choice of $V$.

The geometry of $ \Delta_D $ and more generally the geometry
of the linear rational map
$\PP (( \La \otimes \R)^{\vee} ) \dasharrow \PP ( Herm(m, \C))$  produces
holomorphic invariants of the complex structure.

These are however
related to the topological invariants given
by the linear rational map
$ \PP ( \R^{2d})\cong \PP (( \La \otimes \R)^{\vee} ) \dasharrow \PP (
\La^2(\Ga \otimes
\R)^{\vee} )) \cong
\PP ( \La^2( \R^{2m}))$.
\end{oss}

\section{ Appell Humbert families and Kodaira Spencer map}

We shall consider in this section the family given by the pairs of subspaces
satisfying the Riemann bilinear relations, and the main goal will
be to compute explicitly its tangent space and its Kodaira Spencer map.

\begin{df}
Given $A$ as above, we define $\mathcal
T \mathcal B_A$ as the subset of the product of
Grassmann Manifolds $Gr(m,2m)\times Gr(d,2d)$ defined by

$${\mathcal T}{\mathcal B}_{A}= \{(V,U) \in Gr(m,2m) \times Gr(d,2d)
\ | \ V \cap
\overline{V} = (0),$$
$$U \cap \overline{U} = (0), \ | \ {\rm the \ component \ of}  \ A \ in \
\Lambda^2(\overline{V})^{\vee} \otimes U \ is \ =0 \}.$$
This complex space is called the Appell Humbert space of Torus Bundles.

Since this space is not connected, we restrict ourselves to its
intersection with $\sT_m \times \sT_d$, i.e., we fix respective complex
structures which have the same orientation as  fixed orientations
of $\La$, resp. $\Ga$.
\end{df}

One sees immediately that ${\mathcal T}{\mathcal B}_{A}$ is a complex
analytic variety of codimension at most $dm(m-1)/2$.

Note however that, for $ d \geq 3 , m > > 0$ we get a negative expected
dimension.  The structure of these complex spaces has to be investigated in
general, for our present purposes we limit ourselves to establish the
following

\begin{lem}
If $ d=1$ then the open set ${\mathcal T}{\mathcal B}_{A}
\cap (\sT_m \times \sT_d)$ is connected.
\end{lem}
\Proof

Since  ${\mathcal T}{\mathcal B}_{A}
\cap (\sT_m \times \sT_1)$ fibres onto the upper half plane
$\sT_1$, which is connected, it suffices to show that each fibre is
connected.

Our method of proof will be quite simple minded.
We shall consider the fibration onto $\sT_1$, but without restricting
to $\sT_m$, and considering instead  the larger open
set $\{ V | V \cap
\overline{V} = (0) \}$ (which has two connected components)
and will show that the corresponding fibre has 
exactly two connected components.

Observe then that, if we fix $U$ and write $\La \otimes \C = U \oplus
\overline{U}$,  the connectivity of our space in independent of
the given integral structures, and it amounts to the following problem:
we are given an oriented real vector space $W$ ( $W = \Ga \otimes
\R$) and a real alternating map $ A : \La^2 (W) \ra  U \oplus
\overline{U}$, so that clearly we can write $ A = B \oplus \overline{B}$.

We are seeking for the subspaces
$ V$ such that $\overline{V}$ is isotropic for $B$,
  such that $ W \otimes \C = V \oplus
\overline{V}$, and  such that
  the complex orientation of $V$ induces
the fixed orientation on $W$.

Since the complex linear group is connected,
it is sufficient to show the connectedness of the variety
of frames $ v_1,  \dots v_m $ for $ \overline{V}$.

There are two conditions to be satisfied in the choice
of  $ v_1,  \dots v_m $ :

1) the closed condition that $ v_i$ belongs to the complex
subspace $$ W_i : = \{ v | B (v, v_j) = 0 \ \forall \  1 \leq j \leq i-1 \} ,$$
which has obviously complex dimension at least $ 2m -i + 1$

2) the open condition that $ v_i$ belongs to the open set

$\Omega_i := \{ v_i |  v_1, \overline{v_1}, \dots 
v_i, \overline{v_i}$ are $\C$- linearly
independent $\}$.

We show the connectedness of the variety $\FF_i$ of partial frames
$ v_1,  \dots v_i $ satisfying 1) and 2) by induction 
on $i$.

For $ i=1$, we consider an arbitrary  vector $v_1$,
belonging to the open set 
$\Omega_1 := \{ v_1 |  v_1, \overline{v_1}$ are $\C$- linearly
independent $\}$.

Observe that the complement of $\Omega_1$ has real codimension
$ 2m -1$, thus it does not disconnect a vector space of real dimension
$ 4m$, as soon as $ m \geq 2$.

A similar argument applied to the projection of
$v_i$ into the quotient vector space 
$\tilde{W}_i :=  W \otimes \C / << v_1 , \overline{v_1}, \dots
 v_{i-1} , \overline{v_{i-1}} >>$, which has complex dimension
$ 2m - 2i + 2$, shows that  $\Omega_i$
has a complementary set of real codimension
$ 2m - 2i + 1$, thus we conclude
that 
$\Omega_i \cap W_i$ is connected for $ i < m$.

In the final step instead, we have to remove
the zero locus of the Hermitian form (of $v_m$)
$$ det ( v_1 , \overline{v_1}, \dots
 v_{m} , \overline{v_{m}}) = 0 .$$

We conclude by observing that a Hermitian form has even
positivity and even negativity, thus the complementary
set of its zero locus has exactly two connected components.

Therefore we get, by induction on $i$, that the variety  $\FF_i$ of  
partial frames
$ v_1,  \dots v_i $ satisfying 1) and 2)
is connected for $ i \leq m-1$, while the
variety $\FF_m$ of such frames has at most two connected components
for $ i = m$.

\QED

\begin{df}\label{ah}

The standard (Appell-Humbert) family of torus bundles parametrized by
${\mathcal T}{\mathcal B}_{A}$ is the family of principal holomorphic torus
bundles $X_{V,U}$ over $Y:=Y_V$ and with fibre $T:=T_U$
determined by the cocycle in $H^1(Y, {\mathcal H}(T)_{Y})$ obtained by taking
$f_{\gamma}(v)$ which is the class modulo $\Lambda$ of
$$F_{\gamma}(v):= B'(v, p_V(\ga)) + 2 B''(v, \overline{p_V (\ga)})
    + B''(p_V (\ga), \overline{p_V (\ga)}),\forall v \in V.$$

In other words,  $X_{V,U}$ is the quotient of $  T_U  \times V$
by the action of $\Ga$ such that
$$\ga  ( [u], v)  = (  [u +  B'(v, p_V(\ga)) + 2 B''(v, \overline{p_V (\ga)})
    + B''(p_V (\ga), \overline{p_V (\ga)})], v + p_V (\ga)).$$

\end{df}

\begin{oss}
The above formula differs from the formula given in Definition 6.4 of
\cite{cat04}, where an identification of $\Ga \otimes \R$ with $V$ was
used, and thus $ A(z, \ga)$ was identified with $ B (z, \ga)$.
In the latter formula one had thus

$ -B'(v, \ga) -  B''(v, \overline{p_V (\ga)})$ instead of

    $B'(v, \ga) + 2 B''(v, \overline{p_V (\ga)})
    + B''(p_V (\ga), \overline{p_V (\ga)})$.

However, this alteration does  not affect the proofs of Theorems
6.10,  6.11, 6.12
since in those proofs $B'$
and $B''$ were playing separate roles, and multiplication by $2$ or
$-1$   does not alter
the condition that a certain expression  be zero.

\end{oss}

We also recall from \cite{cat04} the definition of the complete
Appell-Humbert space.

\begin{df}

Given $A$ as above we define

$${\mathcal T}'{\mathcal B}_{A} = \{(V,U, \phi) \ | \ (V,U)
\in {\mathcal T}{\mathcal B}_{A}, \
\phi \in H^1(Y_V, {\mathcal H}(U)_{Y_V}) \cong
\overline{V}^{\vee} \otimes U \}.$$

The complete Appell-Humbert family of torus bundles parametrized by
${\mathcal T}'{\mathcal B}_{A}$ is the family of principal holomorphic torus
    bundles $X_{V,U, \phi}$ on $Y := Y _V$ and with fibre
    $T : = T_U$ determined by the cocycle in $H^1(Y, {\mathcal H}(T)_{Y})$
obtained by taking the sum of $f_{\gamma}(z)$ with the cocycle
$\phi \in H^1(Y_V, {\mathcal H}(U)_{Y_V}) \cong H^1(Y,{\mathcal
O}^d_{Y})$.

\end{df}

Finally we have the following

\begin{teo}{\bf \cite{cat04}}

Any principal holomorphic torus bundle with extension class isomorphic to
$\epsilon \in H^2(\Gamma, \Lambda)$ occurs in the complete Appell-Humbert
family ${\mathcal T}'{\mathcal B}_{A}$.

\end{teo}

Without reproving the above theorem, we want to give a differential geometric
explanation of the cocycle formula above.

The idea is very simple : we consider the homogeneous space
$ M : = \Pi \otimes \R / \Pi$ as a fixed differentiable manifold,
and we use the canonical identification of the tangent space
at the identity in $\Pi \otimes \R $ with
$(\La \otimes \R ) \oplus (\Ga \otimes \R )$ to define
a right invariant  almost complex structure on $ \Pi \otimes \R $ by
translating $ U \oplus V$.

This induces a right invariant almost complex structure on $M$,
and we prove that if the Riemann bilinear relations hold, then
the complex structure is integrable.
Observe that in general $\Pi \otimes \R $ will not have a complex Lie
group structure, because, even if right multiplication is holomorphic,
the inverse mapping needs not be holomorphic.

This interpretation will be quite useful in order to understand
the Kodaira Spencer map of the Appell Humbert families.

\begin{prop}
Let us consider the unique subbundle $T^{(1,0)} \subset T_M \otimes \C$ of the
complexified tangent bundle of
$ M : = \Pi \otimes \R / \Pi$,  which is invariant by right translations, and
is such that
under the identification of the
    tangent space at the identity in $\Pi \otimes \R $ with
$(\La \otimes \R ) \oplus (\Ga \otimes \R )$, we have
$T^{(1,0)}_{Id} = U \oplus V$.

Then, using the other identification of $T_M \otimes \C$ with
the trivial bundle $(\La \otimes \C ) \oplus (\Ga \otimes \C )$,
provided by the  diffeomorphism of $ \Pi \otimes \R$
with $(\La \otimes \R ) \oplus (\Ga \otimes \R )$, with coordinates
$(y,x)$, then

$$T^{(1,0)}_{(y,x)} =  \{(u + A(v, x), v) \ | \ u \in
U, v \in V\} =$$
$$ = \{ u + v + \overline{B''}(v, \overline{p_V(x)}) \ |
    \ u \in U, \ v \in V \} \subset U \oplus
V \oplus \overline{U}$$.

\end{prop}

\Proof
$\forall g = (y,x) \in \Pi \otimes \R$, since $$(y', x') (y,x) =
(y' + y + A (x',x) + 2 S (x,x), x' + x),$$ we have
$$T^{(1,0)}_g=  \{ {R_g}_*(u, v)  \ | \ u \in U, v \in V \} = \{(u +
A(v, x), v) \ | \ u \in
U, v \in V\}
$$

which equals, as claimed,
$$= \{(u + B'(v, p_V(x)) + B''(v, \overline{p_V(x)}) +
\overline{B''}(v,  \overline{p_V(x)}), v)=$$
$$= (u' + \overline{B''}(v,  \overline{p_V(x)}), v)  \ | u' \in U, v \in V\}.$$

\QED

It is now clear that this subbundle $T^{(1,0)} $ is stable under the operation
of the discrete group $\Pi$ given by multiplication from the right.

\begin{REM}
We observe that if we denote by $i: \Pi \otimes \R \ra \Pi \otimes
\R$ the inverse map
sending $g$ to $g^{-1}$ and by $D(i)_g $ the differential of $i$ at the point
$g \in \Pi \otimes \R$, we have  in general that
$D(i)_g ({T_g^{(1,0)}}) \neq
T_{g^{-1}}^{(1,0)}$.

\end{REM}

\Proof
If $g = (y, x) \in \Pi \otimes \R$, $ g^{-1} = i(y,x) = (-y -2S(x,x),
-x)$, so we have
$$D(i)_g(u + A(v, x), v) = (-u -A(v,x) - 4S(x, v), -v).$$
On the other hand
$$T_{g^{-1}}^{(1,0)}  = \{(u + A(v, -x), v) \ | \ u \in U, \ v \in V\}$$
$$\neq \{(u + A(v, x) + 4 S(x,v), v) \ | \ u \in U, \ v \in V\}.$$

\QED

\begin{lem}
\label{integrable}
The above almost complex structure on the differentiable manifold $ M
: = \Pi \otimes
\R / \Pi$ that we have defined is integrable iff the Riemann bilinear relations
are satisfied for the pair $ (U,V)$.

\end{lem}
\Proof

We have already observed the necessity of the Riemann bilinear relations,
so let us just prove the "if" part.

We have $\forall g =(y,x) \in \Pi \otimes \R$, the holomorphic
tangent subbundle
$$T_g^{(1,0)} = \{ u + v +
\overline{B''}(v,
\overline{p_V(x)}) \ | \ u \in U, \ v \in V \} \subset U \oplus V
\oplus \overline{U}$$
and  its conjugate
$$T_g^{(0,1)} = \overline{T_g^{(1,0)}} = \{ \bar{u} + \bar{v} +
B''(\bar{v}, p_V(x)) \ | \ \bar{u} \in \overline{U}, \ \bar{v}  \in
\overline{V } \} \subset
\overline{U} \oplus \overline{V} \oplus U.$$

We have to show that the holomorphic cotangent bundle, which  is the
annihilator of
$T_g^{(0,1)}$, is generated by the differentials of certain local functions.
But we are lucky and we can indeed produce these functions globally.

Consider in fact the diffeomorphism $\Psi: \Pi \otimes \R \ra U \oplus V$
given  by
$$(y,x) \mapsto (p_U(y) - B''(x, p_V(x)),\  p_V(x)) = (p_U(y) -
B''(\overline{p_V(x)},
p_V(x)),\  p_V(x)).$$
$\Psi$ is a diffeomorphism since one can in fact compute the inverse map
explicitly as
$$(u,v) \mapsto (u + \bar{u} + B''(\bar{v}, v) +
\overline{B''(\bar{v}, v)}, v + \bar{v})$$
It is easy to verify that the differentials of the coordinates of
$\Psi$ annihilate
$T_g^{(0,1)}$.

\QED

\begin{cor}

Given $U, V$ satisying the Riemann bilinear relations, the corresponding
integrable complex structure induces a holomorphic action of $\Pi$
on $ U \oplus V$
which corresponds to the cocycle given in definition \ref{ah}.
\end{cor}

\Proof
We use the diffeomorphism $\Psi$ used previously to induce an action
of $\Pi$  on $ U \oplus V$. This action is holomorphic since the complex
structure is invariant by right translation.

We just have to read out
    $ \Psi  \circ R_g \circ  \Psi^{-1} (u,v)$ using the shorthand notation
$g = (y,x)$ and $p_U(y) : = y_U$ , $p_V(x) : = x_V$ (thus, for instance
$ x = x_V + \overline{x_V}$).
We have:
$$ \Psi  \circ R_g \circ  \Psi^{-1} (u,v) =  \Psi  \circ R_g (u +
\bar{u} + B''(\bar{v}, v) +
\overline{B''(\bar{v}, v)}, v + \bar{v})=$$
$$= \Psi  (y + u + \bar{u} + B''(\bar{v}, v) +
\overline{B''(\bar{v}, v)} + A (v + \bar{v} , x), x +  v + \bar{v})= $$
$$= (y_U + u  + B''(\bar{v}, v)
    + B (v + \bar{v} , x)  - B''(x +  v + \bar{v}, x_V +  v), x_V +  v )= $$
$$= (y_U + u  + B' (v , x_V)  + 2 B''( v , \overline{ x_V } )+  B''(
x_V , \overline{ x_V } )   ,
x_V +  v ).
$$

\QED

We give now a description of the tangent space to ${\mathcal
T}{\mathcal B}_{A}$,
observing that   ${\mathcal T}'{\mathcal B}_{A} =
    \{(V,U, \phi) \ | \ (V,U) \in {\mathcal T}{\mathcal B}_{A}, \ \phi
\in \overline{V}^{\vee}
\otimes U\}$ is the restriction to ${\mathcal T}{\mathcal B}_{A}$ of a
    vector bundle on $G(m, 2m) \times G(d, 2d)$.

For $(V,U) \in {\mathcal T}{\mathcal B}_{A}$, the tangent space to
${\mathcal T}{\mathcal B}_{A}$ in $(V,U)$ is a subspace of the tangent space to
the product of the Grassmannians,
$$T_{{\mathcal T}{\mathcal B}_{A},(V,U)} \subset T_{Gr(m,2m) \times Gr(d, 2d),
(V,U)} \cong Hom(V, \overline{V}) \oplus Hom(U, \overline{U}).$$

Let $(L, M) \in  Hom(V, \overline{V}) \oplus Hom(U, \overline{U})$ be
a tangent vector. In order
to determine the Zariski tangent space to
${\mathcal T}{\mathcal B}_{A}$, we work over the ring
of dual numbers $ \C [ \epsilon] / (\epsilon^2)$ and impose the
Riemann bilinear conditions infinitesimally, i.e.,
considering the infinitesimal variations of $V$, resp. $U$, given
     by the graph $V_{\epsilon} $ of
$\epsilon L : V \oplus  \epsilon V \ra  \overline{V}  \oplus
\epsilon \overline{V}$, resp.
the graph $U_{\epsilon}$ of $\epsilon M$.
Concretely,  $V_{\e} = \{ v + \e v' + \e L v | v,v' \in V \}$,  and
    we want  $ A (V_{\e} , V_{\e}   )  \subset U_{\e} $.

$$A(v + \e  v_2+ \e L(v), v' + \e  v'_2+ \e  L(v')) \in U_{\epsilon},
\ \forall v,v' , v_2, v'_2 \in
V.$$ So we must have  elements $ u_0, u'_0 \in U$ such that
$$A(v + \e  v_2+ \e L(v), v' + \e  v'_2+ \e  L(v')) = u_0 + \epsilon
u'_0 + \epsilon Mu_0
,$$

i.e.,
\begin{itemize}
\item
$A(v,v')  = u_0,$
\item
$A(Lv, v') + A(v, Lv') + A(v, v'_2) + A(v_2, v')= u'_0 + Mu_0.$
\end{itemize}

Since $(V,U) \in {\mathcal T}{\mathcal B}_{A}$, $A(v,v') = B(v,v') =
B'(v,v') $,
and the second equation may be rewritten as:
$$B'(v, v'_2) + B'(v_2, v')  -B''(v', Lv) +
\overline{B''(\overline{Lv}, \overline{v'})} +
B''(v, Lv') - \overline{B''(\overline{Lv'}, \overline{v})}
= u'_0 + Mu_0.$$

Finally, observing that  $u'_0 \in U,   Mu_0 \in \overline{U}$, we must have
$$u'_0 =B'(v, v'_2) + B'(v_2, v') -B''(v', Lv) + B''(v, Lv'),$$
$$Mu_0 = MB'(v,v'),$$
$$ MB'(v,v') = \overline{B''(\overline{Lv}, \overline{v'})}   -
\overline{B''(\overline{Lv'}, \overline{v})},$$

and the last one is the equation defining
the tangent space to ${\mathcal T}{\mathcal B}_{A}$ in
$(V,U)$, which we rewrite
using complex conjugation   as

\begin{equation}
\label{tangentspace}
\overline{MB'(v,v')} = B''(\overline{Lv}, \overline{v'})   -
B''(\overline{Lv'},\overline{v}).
\end{equation}

\begin{prop}\label{smooth}

Let $A: \Gamma \times \Gamma \rightarrow \Lambda$ be non zero.

    If $m=2, d=1$, i.e., $\Gamma \cong {\Z}^4$, $\Lambda \cong {\Z}^2$, both
Appell - Humbert spaces
    ${\mathcal T}{\mathcal B}_{A}$ and
${\mathcal T}'{\mathcal B}_{A}$ are smooth.

\end{prop}

{\bf Proof.}

Recall once more that ${\mathcal T}'{\mathcal B}_{A} =
    \{(V,U, \phi) \ | \ (V,U) \in {\mathcal T}{\mathcal B}_{A}, \ \phi
\in \overline{V}^{\vee}
\otimes U\}$ is the restriction to ${\mathcal T}{\mathcal B}_{A}$ of a
    vector bundle on $G(m, 2m) \times G(d, 2d)$ (here $m =2$, $d=1$), so it
suffices to prove that ${\mathcal T}{\mathcal B}_{A}$ is smooth.

We have shown that the tangent space to ${\mathcal T}{\mathcal
B}_{A}$ in $(V,U)$
    is the subspace of $Hom(V, \overline{V}) \oplus  Hom(U, \overline{U})$
    given by kernel of the following linear map

$$G: Hom(V, \overline{V}) \oplus Hom(U, \overline{U}) \rightarrow
Hom(\Lambda^2(V), U),$$

$$(L,M) \mapsto [(v,v') \mapsto  \overline{MB'(v,v')} -
B''(\overline{Lv}, \overline{v'})
    + B''(\overline{Lv'}, \overline{v})],$$
thus it suffices to show that $G$ is surjective.
Since $m =2$ and $d =1$, the dimension of the space
$Hom(\Lambda^2(V), U)$ is one,
    and we only have to prove that the map $G$ is non zero.

For this purpose we observe that if $B' \neq 0$, $G(0,M) \neq 0$, if
$M \neq 0$.

If $B'=0$, we have $B'' \neq 0$ and there must exist $L \in Hom(V,
\overline{V})$
such that $G(L,0) \neq 0$.

In fact, otherwise we would have

$$ B''(\overline{Lv}, \overline{v'})   = B''(\overline{Lv'}, \overline{v}),
\ \forall L \in Hom(V, \overline{V}), \ \forall v,v' \in V.$$

For each $v, v', v'' \in V$, with $v, v'$ linearly independent, there is
      $L$ such that $Lv= \bar{v''} = Lv'$, whence we get
     $B'' (v'', \overline{v' - v}) = 0$ and it follows right away that 
$B'' = 0$.

\hfill{Q.E.D.}

\begin{REM}

The previous proposition contradicts remark 6.13 of \cite{cat04},
which however was due to
a trivial misprint (exchanging $H^1(\hol_X)$ and $H^1(\hol_Y)$).

If $d =1$ and $m \geq 3$, however, ${\mathcal T}{\mathcal B}_{A}$
    is singular at the points where
$B'' =0$.

\end{REM}

{\bf Proof.}

Let $(V, U) \in {\mathcal T}{\mathcal B}_{A}$.

Since $d = dim(U) = 1$, the map

$$G : Hom(V, \overline{V}) \oplus Hom(U, \overline{U}) \rightarrow
Hom(\Lambda^2(V), U),$$

$$(L,M) \mapsto [(v,v') \mapsto  \overline{MB'(v,v')} -
B''(\overline{Lv}, \overline{v'})   + B''(\overline{Lv'},
\overline{v})]$$

has rank one if $B'' =0$.

\hfill{Q.E.D.}

Now we want to determine the Kodaira Spencer map for the Appell-Humbert family.

Recall that $\forall g =(y,x) \in \Pi \otimes \R$, the holomorphic
tangent subbundle is
$$T_g^{(1,0)} = \{ (u + \overline{B''}(v,x), v) \ | \ u \in U, \ v
\in V \} \subset U \oplus V \oplus
\overline{U}$$ thus  its conjugate is
$$T_g^{(0,1)} = \overline{T_g^{(1,0)}} = \{ (\bar{u} + B''(\bar{v},
x),  \bar{v}) \ | \ \bar{u} \in
\overline{U}, \ \bar{v}  \in \overline{V } \} \subset
\overline{U} \oplus \overline{V} \oplus U.$$

Recall that one way to look at the Kodaira Spencer map is to regard it as
associating to $\overline{M} \in
Hom(\overline{U}, U)$,
$ \overline{L} \in Hom(\overline{V}, V)$ a linear map
$ R = \rho (L,M) \in Hom (T_g^{(0,1)},  T_g^{(1,0)} )$ which gives
the infinitesimal
variation of the subspace $T_g^{(0,1)}$ as the  graph of $I + \e  R$, where $I$
denotes the identity
of $T_g^{(0,1)}$.

    So, for $\overline{M} \in
Hom(\overline{U}, U)$,
$ \overline{L} \in Hom(\overline{V}, V)$, we consider the space
$$(T_g^{(0,1)})_{\e }  : = \{ (\bar{u} + \epsilon \overline{M}
\bar{u}+ B''(\bar{v} + \epsilon
\overline{L} \bar{v}, x),
\bar{v} + \epsilon \overline{L} \bar{v} )\ | \ \bar{u} \in
\overline{U}, \ \bar{v}  \in
\overline{V } \}=$$
$$= \{ (\bar{u} + B''(\bar{v}, x) +  \epsilon (\overline{M} \bar{u}+
B''(\overline{L} \bar{v}, x)),
\bar{v} + \epsilon \overline{L} \bar{v} )\ | \ \bar{u} \in
\overline{U}, \ \bar{v}  \in
\overline{V } \}.$$

We obtain $R$ simply by projecting the variation (i.e., the terms
divisible by $\e$)
to $ T_g^{(1,0)}$.

We write this projection using the decomposition:
$ u + \bar{u'} + v + \bar{v'} = $

$= (u -  B''(\bar{v'}, x)) + ( \bar{u'}  -
\overline{B''}(v, x) ) + (v + \overline{B''}(v, x)) + (\bar{v'} +
B''(\bar{v'}, x))$

So we have a map
$$T_g^{(0,1)} \ra  T_g^{(1,0)},$$

and using the isomorphism $U \oplus V \cong T_g^{(1,0)}$ provided by the
choice of $ g = (y,x)$, $ (u,v) \mapsto  u + v +
\overline{B''}(v,x)$, we obtain
$$R : \overline{U} \oplus \overline{V} \ra U \oplus V,$$
$$R : (\bar{u}, \bar{v}) \mapsto (\overline{M} \bar{u}+
B''(\overline{L} \bar{v}, x), \overline{L}
\bar{v} ).$$

We can partly summarize the above description (dependent upon the local choice
of a splitting ) as follows:

\begin{lem}\label{ks}
Consider the exact sequences of vector bundles on $X$:
$$ 0 \ra U  \ra   T^{(1,0)} \ra V \ra 0, \ \  0 \ra \overline{U} \ra
T^{(0,1)} \ra
    \overline{V} \ra 0.$$
Then if $R = \rho (L,M, \phi)$ is the Kodaira Spencer image of the
tangent vector
$ (L,M, \phi)$ to the complete Appell Humbert space, then the image
of $ R$ in $Hom ( T^{(0,1)}, V)$ is the image of $\overline{L} \in
Hom(\overline{V}, V)$,
and the image
of $ R$ in $Hom ( \overline{U}, T^{(1,0)})$ is the image of
$\overline{M} \in Hom(\overline{U}, U)$.
\end{lem}

\section{ On the versality of the complete Appell Humbert family }

We recall for the reader's benefit some results from \cite{cat04}:

\begin{teo}
The cokernel of $ 0 \ra H^0 (\Omega^1_Y )\ra H^0 (\Omega^1_X )$ is the
subspace of ${U}^{\vee}$ which annihilates the image of the Hermitian
part of $A$, i.e.,  of the component $B''$ in
$[( V \otimes \bar{V})^{\vee}  \otimes (U )].$

It follows in particular that $X$ is parallelizable if and only if the
Hermitian part of $A$ is zero.
\end{teo}

\qed
\begin{cor}
The space $H^0 ( d\hol_X)$ of closed holomorphic 1-forms on $X$ contains
the pull- back of $H^0 ( \Omega^1_Y)$ with cokernel the subspace $U^*$
of $U^{\vee}$ which annihilates the image of $B$, i.e.,  $U^* = \{ \beta \ | \
\beta
\circ B (z,\gamma)= 0, \
\forall
\gamma, \ \forall z\} .$
\end{cor}

\begin{teo}
The cokernel of $ 0 \ra H^1 (\hol_Y )\ra H^1 (\hol_X )$ is the subspace
of $\bar{U}^{\vee}$ which annihilates the image of the anti-complex component
    of $A$, i.e., of the conjugate of the component $B'$ in
$[( \Lambda^2 V )^{\vee}  \otimes (U )].$
\end{teo}

In order to analyse the problem of describing the small deformations of
principal holomorphic torus bundles over tori we need to calculate the
cohomology groups of the tangent sheaf $\Theta_X$.

These fit into the following exact sequences, according to

\begin{cor}{\bf \cite{cat04}}
$H^i (\Theta_X)$ fits into a short exact sequence
$$ 0 \ra coker b_{i-1} \ra H^i (\Theta_X) \ra ker b_{i} \ra
0,$$ where  $b_i : V \otimes H^{i} (\hol_X) \ra U \otimes H^{i+1}
(\hol_X)$ is given by cup product and contraction with $ B'' \in  [(  \bar{V})
\otimes V )^{\vee} \otimes (U )]$.
\end{cor}

We want to give now a more explicit description of these groups
for $i=1$.

We have an exact sequence

$$0 \ra U \otimes \hol_X \ra \Theta_X \ra V \otimes \hol_X \ra 0 ,$$
whose extension class is provided by $B''$,  yielding a cohomology
exact sequence
$$(**)  V \ra   U \otimes H^1(X, \hol_X)   \ra H^1(X, \Theta_X ) \ra
V \otimes H^1(X, \hol_X) \ra
U \otimes H^2(X, \hol_X) , $$
and the above cited corollary says that the first and the last maps
are given by cup product and
contraction with $ B'' $.

The following proposition allows us to describe the cohomology groups
$H^1(X, \hol_X)$,
$H^2(X, \hol_X)$.

\begin{prop}{\cite{cat04}}
The Leray spectral sequence for the sheaf $\hol_X$ and for the map
$f: X \ra Y$ is a
spectral sequence  which degenerates at the $E_3$ level and with
$E_2$ term  $ = (H^i
(\RR^j f_*
\hol_X ) , d_2 )$, where $d_2 : (H^i (\RR^j f_* \hol_X )= H^i (
\Lambda^j (\bar{U}^{\vee}) \otimes \hol_Y) =
\Lambda^i (\bar{V}^{\vee}) \otimes \Lambda^j (\bar{U}^{\vee}) \ra \Lambda^{i+2}
(\bar{V}^{\vee}) \otimes \Lambda^{j-1} (\bar{U}^{\vee})$ is provided
by cup product
    and contraction with $\bar{B'} \in \Lambda^2 (\bar{V}^{\vee})
\otimes (\bar{U})$.
\end{prop}

The Leray spectral sequence yields then an exact sequence

$$0 \ra H^1(Y, \hol_Y) = \overline{V}^{\vee} \ra H^1(X, \hol_X) \ra
H^0(\RR^1 f_*
\hol_X )  = \overline{U}^{\vee} \ra $$
$$\ra H^2(Y, \hol_Y) = \Lambda^2(\overline{V}^{\vee}) \ra
H^2(X,
\hol_X),$$ and the map
$$ H^0(\RR^1 f_*\hol_X ) = \overline{U}^{\vee}\ra H^2(Y, \hol_Y) =
\Lambda^2(\overline{V}^{\vee})$$
is given by $\overline{B'}$.

So we  see again that $H^1(X, \hol_X) \cong \overline{V}^{\vee} \oplus
Ker(\overline{B'}), $
where we view $\overline{B'} $ as a map
$\overline{B'}: \overline{U}^{\vee}\ra \Lambda^2(\overline{V}^{\vee})$.

Thus,  to determine $H^1(X, \Theta_X )$, it remains to compute
$H^2(X, \hol_X)$.

Since the Leray spectral sequence degenerates at $E_3$ level, we have

$$H^2(X, \hol_X) \cong E_3^{0,2} \oplus E_3^{2,0} \oplus E_3^{1,1}, $$
where
$E_3^{0,2} = Ker(d_2: E_2^{0,2} \ra E_2^{2,1} ) =
Ker(  \overline{B'} :\Lambda^2(\overline{U}^{\vee}) \ra
\overline{U}^{\vee} \otimes
\Lambda^2(\overline{V}^{\vee}))$,

$E_3^{1,1} = Ker(d_2: E_2^{1,1} \ra E_2^{3,0} ) =
Ker(\overline{B'} :\overline{U}^{\vee} \otimes  \overline{V}^{\vee} \ra
\Lambda^3(\overline{V}^{\vee}))$,

$E_3^{2,0} = \frac{E_2^{2,0}}{Im(d_2: E_2^{0,1} \ra E_2^{2,0})}
= Coker (\overline{B'} :  \overline{U}^{\vee}\ra
\Lambda^2(\overline{V}^{\vee}))$.

Now we want to determine the kernel of the coboundary map
$$F: V \otimes H^1(X, \hol_X) \ra U \otimes  H^2(X, \hol_X),$$
in the exact sequence $(**)$

$$V \ra U \otimes H^1(X, \hol_X) \ra H^1(X, \Theta_X) \ra V \otimes
H^1(X, \hol_X) \ra U \otimes  H^2(X, \hol_X).$$

in terms of the isomorphism
$$V \otimes H^1(X, \hol_X) \cong V \otimes (\overline{V}^{\vee}
\oplus Ker(\overline{B'})) =
    (V \otimes \overline{V}^{\vee})
\oplus  (V \otimes Ker(\overline{B'})),$$

\begin{lem}
The kernel of the coboundary map
$$F: V \otimes H^1(X, \hol_X) \ra U \otimes  H^2(X, \hol_X),$$
in the exact sequence $(**)$ consists of the pairs
$(\bar{L}, \sum_i v_i \otimes  \overline{u_i}^{\vee})$ such that

1)  there exists $\overline{M} \in
\overline{U}^{\vee} \otimes U$ such that
$B''(\overline{L}) =  \overline{M} \circ \overline{B'} :
\La^2(\overline{V})
\ra  U $,

2)  $\sum_i B'' (v_i )\otimes  \overline{u_i}^{\vee} = 0$ in
$ U \otimes \overline{V}^{\vee}
\otimes\overline{U}^{\vee} $ for each choice of
$v_i \in V$,
$\overline{u_i}^{\vee} \in  Ker(\overline{B'})$
\end{lem}

\Proof

Since $F$ is given by cup product and contraction with $B''$, for
    $\overline{L} \in V \otimes \overline{V}^{\vee}$, $v \in V$,
$\overline{u}^{\vee} \in  Ker(\overline{B'})$,  we have

$$F(\overline{L}, v \otimes \overline{u}^{\vee}) =
(\zeta(\overline{L}), B''(v)
    \otimes \overline{u}^{\vee}),$$
where
$\zeta(\overline{L}) = [B''(\overline{L})] \in \frac{\La^2(\overline{V}^{\vee})
\otimes U}{ Im(\overline{B'}) \otimes U} = E_3^{2,0} \otimes U,$

$( B''(v)  \otimes \overline{u}^{\vee}) \in U \otimes
    Ker(\overline{B'}_{|  \overline{U}^{\vee}
\otimes\overline{V}^{\vee} })  = U \otimes E_3^{1,1}$.

So $ \zeta(\overline{L}) = 0$ if and only if
$B''(\overline{L}) \in Im(\overline{B'} \otimes Id:
\overline{U}^{\vee} \otimes U
    \ra \La^2(\overline{V}^{\vee}) \otimes U)$, i.e.
    there exists $\overline{M} \in
\overline{U}^{\vee} \otimes U$ such that
$B''(\overline{L}) = \overline{M} \circ \overline{B'}  :
\La^2(\overline{V})
\ra \overline{U} \ra U $.

\QED

\begin{oss}
Observe that condition 1)  says exactly that the pair $(L, M)$
belongs to the tangent space
    to the Appell Humbert space.
    \end{oss}

\begin{prop}\label{KS}
The Kodaira Spencer map of the complete Appell Humbert family is surjective
under one of the following assumptions:

i] $h^1(\hol_X) = h^1(\hol_Y) $, i.e., $Ker(\overline{B'}) = 0$

ii] $d=1$ and $B' \neq 0$

iii] $d=1$ and  $B (v, x)$ non degenerate in the first factor.

iv] $X$ is parallelizable and $A$ is surjective
as a linear map between $\La^2 (V \oplus \overline{V}) \ra ( U \oplus
\overline{U})$.
\end{prop}
\Proof

By the previous remark and by lemma \ref{ks} the image  of $\rho (L,M, \phi)$
inside $H^1 (V \otimes \hol_X)$ is simply the pair $(\bar{L}, 0)$.

Thus, we first observe that the subspace of the tangent vectors
with $L=0$, namely, $\{ (0, M, \phi ) | M \circ B' = 0  \}$ maps
onto $ \{ (\overline{M}, \phi ) | \overline{M} \circ \overline{B'} = 0  \}
=  (U \otimes  ker
(\overline{B'})) \oplus ( U \otimes
\overline{V}^{\vee}) = H^1 (U \otimes \hol_X)$.

Second, by condition 1) above, the image of the Kodaira Spencer map of
the complete Appell
Humbert family is surjective onto $ker (F)  \subset H^1 (V \otimes
\hol_X)$ if and
only if $ker (F) \subset (V \otimes \overline{V}^{\vee}) $.

Therefore, we want that $\sum_i B'' (v_i )\otimes
\overline{u_i}^{\vee} = 0$ in
$ U \otimes \overline{V}^{\vee}
\otimes\overline{U}^{\vee} $ with
$v_i \in V$,
$\overline{u_i}^{\vee} \in  Ker(\overline{B'})$ implies $\sum_i v_i \otimes
\overline{u_i}^{\vee} = 0$.

Therefore this holds trivially under condition i],
that $Ker(\overline{B'}) = 0$.

In turn, condition ii] : $d=1$ and $B' \neq 0$ implies condition i].

In the case  $d=1$   we do not need to take a sum of simple tensors,
but just one,
and then $ B'' (v )\otimes
\overline{u}^{\vee} = 0$, for $\overline{u}^{\vee} \neq 0$, holds iff
$ B'' (v )=0 \in
( \overline{V}^{\vee} \otimes U) $.
However, since we may assume $B'=0$,  the last condition  is
equivalent to $ B (v) = 0 $,
which contradicts iii].

Finally, we show that iv] implies i]: in fact we know that parallelizability
is equivalent to the vanishing of $B''$. Whence, we have $ A = B'
\oplus \overline{B'}$,
thus $A$ is surjective if and only if $B'$ is surjective as a linear map
$\La^2 (V ) \ra ( U )$. But then the conjugate of the transpose
is injective, i.e., condition i] holds.

\QED

\begin{oss}

On a compact complex manifold $X$ one has the Kodaira
inequalities
$$2 h^1 (X,  \hol_X) \geq h^1 (X,  \hol_X) + h^0( d {\hol_X })
\geq b_1(X),$$
which hold iff $X$ has a very good (cf. e.g. \cite{cat04}) Albanese
variety, of dimension
$h^1 (X,  \hol_X) = h^0( d {\hol_X }) = \frac{1}{2}
    \ b_1(X)$.

For principal holomorphic torus bundles then $h^1(\hol_X) = h^1(\hol_Y)$
implies by semicontinuity that every small deformation $X_t$ deforms
together with an Albanese map $X_t \ra Y_t$. If $d=1$  necessarily
$X_t \ra Y_t$ is an elliptic bundle, but for $d \geq 2$ we only have
a differentiable torus bundle as a deformation, and therefore
i] of the previous proposition yields a further contribution.

\end{oss}

Recall the following
\begin{teo}{\bf \cite{cat04}}
Let $ f:  X \ra Y$ be a holomorphic principal bundle with base a complex torus
$ Y _V$ of dimension $m$, and fibre an elliptic curve $T_U$.

Assume moreover $ \pi_1 (X) := \Pi $ to be a nontrivial  central extension
$$ 1 \ra \Lambda \ra \Pi \ra \Gamma \ra 1$$
classified
by  a cohomology class
$\epsilon  \neq 0 \ \in H^2(Y, \Lambda )$
whose associated bilinear form $A$ has an image of dimension $=2$.

Then every limit of  manifolds in the complete Appell-Humbert family is again a
holomorphic principal bundle $ f' :  X' \ra Y'$ with fibre an
elliptic curve $T'$,
and thus occurs in the complete Appell-Humbert family .

\end{teo}

The above result combines with the following one:

\begin{teo}
    Let $ f:  X \ra Y$ be a holomorphic principal bundle with base a
complex torus
$ Y: = Y _V$ of dimension $m = 2$, and fibre an elliptic curve $T : = T_U$.

Then the complete Appell Humbert family is versal at $X$ if
either $B' \neq 0$, or $B'=0$ and  $B (v, x)$ non degenerate in the
first factor,
or $X$ is parallelizable and $A \neq 0$.
\end{teo}
\Proof
By Proposition \ref{smooth} the complete Appell Humbert space is smooth
for $m=2,d=1$, thus versality would follow from the surjectivity of
the Kodaira Spencer map.

Surjectivity under the first two respective assumptions follows directly from
ii] and iii] of proposition \ref{KS}, whereas, if $X$ is parallelizable,
we argue as in iv] of \ref{KS}: $B'' = 0$, and if also $B' = 0$, then $A = 0$,
a contradiction. Thus $B' \neq 0$ and we are done.

\QED

\begin{prop}
Assume $m=2$, $d=1$ and assume further that $A$ is non degenerate.
Then  the case : there is $(U,V)$ such that $B'=0$
and $B (v, x)$ is degenerate
in the first factor occurs only if the associated pencil

$ \PP^1 \cong \PP ( (\La \otimes \R)^{\vee} ) \into \PP ( \La^2(\Ga \otimes
\R)^{\vee}) \cong
\PP ( \La^2( \R^{4}))$
intersects the Klein Pfaffian quadric in two distinct complex
conjugate points.
\end{prop}

\Proof
Assume that $ v \in V$ is such that $ B(v,x) = 0$ for each $x$, and
recall the assumption  $ B' = 0$. We can choose $w \in V$ such
that $ \{ v, w \}$ is a basis of $V$.

We have $B'' ( v, \bar{v}) = 0 = \overline{ B''  }( v, \bar{v})$,
$B'' ( v, \bar{w}) = 0 = \overline{ B''  }( w, \bar{v})$.

We  set $  B''  ( w, \bar{v}): = \zeta$ so that
$\overline{ B''  }( v, \bar{w}) = - \bar{\zeta}$,
and observe that  $\zeta \neq 0$,
otherwise $v$ is in the kernel of $A$.

Up to a change of basis replacing $w$ by $ w + \mu v$,
we may assume that also $B'' ( w, \bar{w}) = 0 = \overline{ B''  }(
w, \bar{w})$.

We have found two real isotropic subspaces for $A$, whose complexifications
are spanned by
$\{ v, \bar{v}\}$,$ \{ w, \bar{w} \}$ respectively.

These two subspaces are spanned by
$\{ v \pm \bar{v}\}$,$ \{ w \pm \bar{w} \}$ and we compute the corresponding
values of $B''$ on pairs of such vectors.
Setting $ \e_1, \e_2 = \pm 1$ we get
$$ B''(  v + \e_1 \bar{v},   w  + \e_2 \bar{w}) =  - \e_1 \ze$$
$$ \overline{ B''  }(  v + \e_1 \bar{v},   w  + \e_2 \bar{w}) =  -
\e_2 \bar{\ze}.$$

We observe at this point that $ e_1 :=  \ze + \bar{\ze}$, $e_2 :=
\frac{1}{i}( \ze - \bar{\ze})$
give a real basis for $\La \otimes \R$,
whence in the basis $\{ 2 Re(v), 2 Im (v), 2 Re (w), 2 Im (w) \}$
$A$ is given by a $ 2 \times 2$  block matrix  with diagonal
blocks equal to zero and with upper right block equal to:
$$M=\left(
\begin{array}{cc}
- e_1& - e_2\\
e_2 & -e_1
\end{array}
\right).
$$
Therefore, the corresponding Pfaffian of $\mu_1 A_1 + \mu_2 A_2$,
if $A_1, A_2$ are the components of $A$ with respect to the basis
$e_1, e_2$, equals $\mu_1^2 + \mu_2^2$, and it has no real roots.

\QED

Combining the last three stated results follows rightaway the following

\begin{teo}
Consider the family of holomorphic principal torus bundles $ f:  X \ra Y$
with  base a complex torus
$ Y _V$ of dimension $2$, and with  fibre an elliptic curve $T_U$,
corresponding to a bilinear form $A$ such that

1) $A$ is non degenerate and ${\rm Im}  A$ has dimension $2$

2) the associated pencil

$ \PP^1 \cong \PP ( ( \La \otimes \R)^{\vee} ) \into \PP
   ( \La^2(\Ga \otimes
\R)^{\vee}) \cong
\PP ( \La^2( \R^{4}))$
intersects the Klein Pfaffian quadric in at least one real point.

Then this family forms a connected component of the  Teichm\"uller space.

\end{teo}

\begin{oss}
One example of an elliptic bundle over
a 2-dimensional torus is the so called {\bf Iwasawa manifold} (cf.
\cite{k-m71}).
It is parallelizable, being the quotient of the complex Lie group
    $N$ of $3 \times 3$ upper trangular
matrices with all eigenvalues equal to $1$ ($N \cong \C^3$  as a
complex manifold)
by the discrete cocompact subgroup $\Pi$ which is the subgroup of the
matrices with entries in the subring $\Lambda \subset \C$,  $\Lambda : =
\Z [i]$.

Nakamura observed that the Kuranishi family of $X$ is smooth of dimension $6$,
and there are small
deformations which are not parallelizable. This follows by our result in
greater generality.

In this case, the alternating form $A$ is gotten as the
antisymmetrization of the product map $\C \times \C \ra \C$ , which
induces  $\Lambda \times \Lambda \ra \Lambda$.

One can explicitly calculate the integral bilinear map $A$, and
it turns out that there are indeed cases where $B'=0$ and
$ B(v,x)$ is degenerate in the first factor.

In fact, one can verify that the corresponding  pencil of
alternating forms intersects the Pfaffian Klein quadric in two
distinct complex conjugate points. Thus, there are
manifolds where the Kodaira Spencer image of the
complete Appell Humbert family is not surjective.

If the complete Appell Humbert family does not yield
an open set in the Kuranishi space, then the following should happen:
we have a $X$ for which $h^1(\hol_X) = 3, h^0 (d \hol_X) = 2=
h^0(\Omega^1_X)$. But there is a small deformation
$X_t$ such that   $h^1(\hol_{X_t}) = 3, h^0 (d \hol_{X_t}) = 1$.

It is interesting to see whether this can happen,
for instance if this could be related to the existence
of a fibration of  $X$  over an elliptic
curve which deforms along, producing
another component of the Kuranishi space.

\end{oss}

We conclude by  observing, as recalled to us by Marian Aprodu:

\begin{oss}
    Let $ f:  X \ra Y$ be a holomorphic principal bundle
with base a complex torus
$ Y: = Y _V$ of dimension $m $, and fibre a torus $T : = T_U$, such that
$A \neq 0$.

Then $X$ is not K\"ahler.

\end{oss}

Indeed, this was proven by Blanchard in \cite{bla56}
more generally for holomorphic fibre
bundles whose transgression $ H^1 ( T , \R) \ra H^2 ( Y , \R)$
is nonzero ( cf. also \cite{hoefer}).

\bigskip

{\bf Acknowledgements.}

The  authors wish to thank Marian Aprodu, resp. Alessandro Ghigi,
    for a useful conversation.

\vfill

\noindent
{\bf Author's address:}

\bigskip

\noindent
Prof. Fabrizio Catanese\\
Lehrstuhl Mathematik VIII\\
Universit\"at Bayreuth\\
    D-95440, BAYREUTH, Germany

e-mail: Fabrizio.Catanese@uni-bayreuth.de

\noindent
Dr. Paola Frediani\\
Dipartimento di Matematica\\
Universit\`a di Pavia\\
    I-27100 Pavia, Italy

e-mail: paola.frediani@unipv.it


\bigskip



\begin{thebibliography}{99}

\bibitem[Apr]{apr}
M. Aprodu,
``{\em An Appell-Humbert theorem for hyperelliptic surfaces}'',
{\bf J. Math. Kyoto Univ. 38, no. 1} (1998), 101-121.

\bibitem[BPV84]{bpv84}
W. Barth, C. Peters and A. Van de Ven,
``{\em Compact Complex Surfaces}'', {\bf Ergebnisse der Mathematik und
ihrer Grenzgebiete, 3.F, B. 4, Springer-Verlag}, 1984.


\bibitem[Bla53]{bla53}
A. Blanchard, ``{\em Recherche de structures analytiques complexes
sur certaines vari\'et\'es}'', {\bf C.R. Acad.Sci., Paris 238} (1953),
657-659.

\bibitem[Bla56]{bla56}
A. Blanchard, ``{\em Sur les vari\'et\'es
analytiques complexes}'', {\bf Ann.Sci. Ec. Norm. Super., III Ser., 73}
(1956), 157-202.


\bibitem[Cal58]{cal58}
E. Calabi,
``{\em Construction and properties of some 6-dimensional almost complex
manifolds}, {\bf Trans. Amer. Math.Soc. 87  } (1958), 407-438.





\bibitem[Cat91]{cat91}
F. Catanese,
``{\em Moduli and classification of irregular K\"ahler manifolds
    (and algebraic varieties) with Albanese general type fibrations.
    Appendix by Arnaud Beauville.}''
{\bf Inv. Math. 104} (1991) 263-289; Appendix 289 .


\bibitem[Cat95]{cat95}
F. Catanese,
``{\em Compact
complex manifolds bimeromorphic to tori.}'' Proc. of the Conf. "Abelian
Varieties", Egloffstein
1993, {\bf  De Gruyter }(1995), 55-62.





\bibitem[Cat02]{cat02}
F. Catanese,
``{\em Deformation types of  real and complex manifolds}'', Proc. of
the Chen-Chow
Memorial Conference, Contemporary trends in algebraic geometry
and algebraic topology (Tianjin,
2000),
      {\bf Nankai Tracts Math., 5, World Sci. Publishing, River Edge, NJ}
      (2002), 195--238 .

\bibitem[Cat04]{cat04}
F. Catanese,
``{\em Deformation in the large of some complex manifolds, I}'',
   {\bf Ann. Mat. Pura Appl. (4)  183, no. 3 }, Volume in Memory of
Fabio Bardelli, (2004),   261--289 .



\bibitem[Che58]{che58}
S.S. Chern,
``{\em Complex manifolds}'',
{\bf Publ. Mat. Univ.  Recife} (1958).


\bibitem[FG65]{fg65}
W. Fischer, H. Grauert,
``{\em Lokal-triviale Familien kompakter komplexer Mannigfaltigkeiten }''
{ \bf Nachr. Akad. Wiss. G\"ottingen, II. Math.-Phys. Kl. 1965} (1965), 89-94 .



\bibitem [Hir62]{hir62}
H. Hironaka,
``{\em An example of a non-K\"ahlerian complex-analytic deformation of
K\"ahlerian complex structures}'', {\bf Ann. Math. (2) 75,}
(1962),190-208.

\bibitem[H\"o93]{hoefer}
T. H\"ofer,
   ``{\em  Remarks on torus principal bundles}'', {\bf  J. Math. Kyoto Univ.
33 ,  no. 1} (1993), 227--259.



\bibitem [Ko64]{ko64}
K. Kodaira ,''{\em On the structure of complex analytic
surfaces I}'', {\bf Amer. J. Math.  86} (1964), 751-798 .

\bibitem [Ko68]{ko68}
K. Kodaira ,''{\em On the structure of complex analytic
surfaces IV}'', {\bf Amer. J. Math.  90} (1968), 1048-1066 .

\bibitem [K-M71]{k-m71}
K. Kodaira, J. Morrow ,``{\em Complex manifolds}'' {\bf Holt,
Rinehart and Winston},  New York-Montreal, Que.-London (1971).

\bibitem[K-S58]{k-s58}
K. Kodaira , D. Spencer ''{\em On deformations of complex analytic
structures I-II}'', {\bf Ann. of Math.  67} (1958), 328-466 .

    \bibitem[Ku1]{ku1}
M. Kuranishi,
    ''{\em On a type of family of complex structures}'', {\bf Ann.
of Math. (2) 74 } (1961), 262--328.

    \bibitem[Ku2]{ku2}
M. Kuranishi,
    ''{\em New proof for the existence of locally complete families of complex
structures}'',  Proc. Conf. Complex Analysis (Minneapolis, 1964)
{\bf Springer, Berlin} (1965), 142--154.

\bibitem[Nak75]{nak75}
I. Nakamura,
''{\em Complex parallelisable manifolds and their small deformations}''
{ \bf J. Differ. Geom. 10}, (1975),85-112 .

\bibitem[Nak98]{nak98}
I. Nakamura,
''{\em Global deformations of $\PP^2$-bundles over $\PP^1$}''
{ \bf J. Math. Kyoto Univ. 38, No.1} (1998), 29-54 .




\bibitem[Somm75]{somm75}
A. J. Sommese ,``{\em Quaternionic manifolds}'',
{\bf
Math. Ann. 212}  (1975), 191-214.


     \bibitem [Su70]{su70}
     T. Suwa, ``{\em On hyperelliptic surfaces}'', {\bf J. Fac. Sci.
Univ. Tokyo, 16} (1970), 469-476 .


    \bibitem [Su75]{su75}
     T. Suwa, ``{\em Compact quotient spaces of $\C^2$ by affine transformation
groups}'', {\bf J. Differ. Geom. 10},  (1975) 239-252.




\bibitem [Su77-I]{su77-I}
     T. Suwa, ``{\em Compact quotients of $\C^3$ by affine transformation
groups, I.}'', {\bf Several complex Variables, Proc. Symp. Pure Math.
30, Part 1,
Williamstown 1975,} (1977), 293-295.




\bibitem [Su77-II]{su77-II}
     T. Suwa, ``{\em Compact quotients of $\C^3$ by affine transformation
groups, II.}'', {\bf Complex Analysis and Algebraic Geometry, Iwanami
Shoten} (1977), 259-278 .

\bibitem [Ue75]{ue75}
     K. Ueno, ``{\em
Classification theory of algebraic varieties and compact complex spaces.}

{\bf Lecture Notes in Mathematics 439,}
Springer-Verlag, XIX, 278 p.(1975).

\bibitem [Ue80]{ue80}
     K. Ueno, ``{\em On three-dimensional compact complex manifolds with non
positive Kodaira dimension}'',  {\bf Proc. Japan Acad.,
vol.56,S.A.n.10} (1980),
479-483.

    \bibitem [Ue82]{ue82}
     K. Ueno, ``{\em Bimeromorphic Geometry of algebraic and analytic
threefolds}'', in 'C.I.M.E. Algebraic Threefolds, Varenna 1981'
{\bf Lecture
Notes in math. 947} (1982), 1-34 .

\bibitem [Ue87]{ue87}
     K. Ueno, ``{\em On compact analytic Threefolds with non trivial Albanese
Tori}'',  {\bf Math. Ann. 278} (1987), 41-70.




\end{thebibliography}
\end{document}